\documentclass[sn-mathphys-num]{sn-jnl}

\usepackage{graphicx}%
\usepackage{multirow}%
\usepackage{amsmath,amssymb,amsfonts}%
\usepackage{amsthm}%
\usepackage{mathrsfs}%
\usepackage[title]{appendix}%
\usepackage{xcolor}%
\usepackage{textcomp}%
\usepackage{manyfoot}%
\usepackage{booktabs}%
\usepackage{algorithm}%
\usepackage{algorithmicx}%
\usepackage{algpseudocode}%
\usepackage{listings}%
\usepackage{dsfont}

\usepackage{tikz}      

\usepackage{pgfplots}
\pgfplotsset{width=6.5cm,compat=1.9}

\usepackage{hyperref}   
\hypersetup{
    colorlinks=true,
    linkcolor=purple,
    filecolor=magenta,      
    urlcolor=cyan,
    citecolor=teal
    }

\newcommand{\mZ}{\mathcal{Z}}
\newcommand{\mD}{\mathcal{D}}

\newcommand{\mF}{\mathcal{F}}

\newcommand{\mC}{\mathcal{C}}

\newcommand{\ux}{\underline{x}}

\newcommand{\uh}{\underline{h}}

\newcommand{\uv}{\underline{v}}
\newcommand{\uz}{\underline{z}}

\newcommand{\N}{\mathbb{N}}
\newcommand{\Sf}{\mathbb{S}}
 
\newcommand{\T}{\mathbb{T}}

\newcommand{\R}{\mathbb{R}} 

\newcommand{\Lp}{\mathbb{L}}

\newcommand{\ind}{\mathds{1}}

\newcommand{\dd}{\mathrm{d}}

\newcommand{\lnorm}{\left\|}
\newcommand{\rnorm}{\right\|}

\newcommand{\CQFD}{\hfill $\square$}

\newcommand{\und}{\mbox{ and }}

\theoremstyle{plain}

\newtheorem{theo}{Theorem}[section]
\newtheorem{prop}{Proposition}[subsection]

\theoremstyle{remark}
\newtheorem*{prof}{Proof}
\newtheorem*{rema}{Remark}

\newcommand{\supess}{\mathop{\smash{\mathrm{supess}}}}

\everymath{\displaystyle} 

\usepackage{anyfontsize}

\begin{document}

\title[Article Title]{On the derivation of the linear Boltzmann equation from the nonideal Rayleigh gas}

\subtitle{Adaptative pruning and improvement of the convergence rate}


\author[1]{\fnm{Florent} \sur{Foug\`{e}res}}\email{florent.fougeres@ens.fr}

\date{2024}

\affil[1]{\orgdiv{DMA}, \orgname{\'{E}cole normale sup\'{e}rieure}, \orgaddress{\street{45 rue d'Ulm}, \city{Paris}, \postcode{75005},  \country{France}}}

\abstract{This paper's objective is to improve the existing proof of the derivation of the Rayleigh--Boltzmann equation from the nonideal Rayleigh gas~\cite{2016brownian}, yielding a far faster convergence rate. This equation is a linear version of the Boltzmann equation, describing the behavior of a small fraction of tagged particles having been perturbed from thermodynamic equilibrium. This linear equation, derived from the microscopic Newton laws as suggested by the Hilbert's sixth problem, is much better understood than the quadratic Boltzmann equation, and even enable results on long time scales for the kinetic description of gas dynamics.

The present paper improves the physically poor convergence rate that had been previously proved, into a much more satisfactory rate which is more than exponentially better.}


\ackno{I send my many thanks through these lines to my PhD directors Isabelle Gallagher and Sergio Simonella, for their tireless help and support, as well as to all the persons who partake in making the DMA a living, stimulating and cosy place.}



\maketitle  

\newpage

\section{Introduction}

The Boltzmann equation, introduced in 1872 by Ludwig Eduard Boltzmann~\cite{1872boltzmann}, provides a model for rarefied gas dynamics that has paved the way to a flourishing litterature exploring the corresponding mesoscopic scale. In fact, this kinetic representation may be used as an intermediate step in the derivation of fluid mechanics equations~(see among many others~\cite{1991BGL, 1993BGL,2001ns,2008euler}) and has permitted various numerical applications~\cite{2022lbm, 2023lbm}, along with a theoretical comprehending of  the intrinsic behavior of such rarefied gases. Indeed, very soon after the equation's formulation, introducing the concept of entropy, Boltzmann has shown that the solutions to this equation irreversibly converge, for long times, towards an equilibrium. This equilibrium is well known: it is uniform in space and distributed according to the Maxwell Gaussian distribution~(\ref{def:max}) in velocities, which depends only on the energy, or temperature, of the system.

On the other hand, the microscopic state of the gas, from which the Boltzmann equation is derived in the very low density limit, is given by classical Newton equations, and the solutions to such equations are completely time reversible; this seeming paradox naturally led Boltzmann's contemporaries into doubting the validity of his model.
Nevertheless, the rigorous derivation of this mesoscopic Boltzmann equation from microscopic Newton equations has eventually been proved mathematically  in 1975 by Oscar Erasmus Lanford III~\cite{1975lanford, 1976lanford} in the case of hard sphere interactions; but the methods he used suffer from a strong ridigity that hinders to extend his result for long time scales. In fact, his proof is only valid for very small times, when only about a fifth of particles have collided. The major obstruction is the correlation that happens when two particles collide: the system's chaotic properties deteriorate over time, making it very hard to deal with any recollision of particles.

Nonetheless, close to thermodynamic equilibrium, the statistical stability of the dynamics guarantees the propagation of a certain amount of chaos and provides a control of correlations in a very strong sense. Hence, in the Rayleigh~gas model, describing the behavior of a small fraction of tagged particles near equilibrium~\cite{1988spohnkineq}, a very similar proof leads to the derivation \emph{for long time scales} of a linear version of the Boltzmann equation, called Rayleigh--Boltzmann equation; this has been the work of Henk van Beijeren, Lanford, Joel Louis Lebowitz and Herbert Spohn in 1980~\cite{1980linear}.

A few decades later, in 2013, Isabelle Gallagher, Laure Saint-Raymond and Benjamin Texier reopened the work of  Lanford and his former student, Francis Gordon King~\cite{1975king}, filling some gaps including the case of compactly supported potentials, and hence providing precise estimates on the convergence rate in Lanford's theorem. Eventually, this work has led in 2016 to an article by Thierry Bodineau, Gallagher and Saint-Raymond on the convergence rate in the linear case, and its dependence on the long time scaling, so as to infer Brownian hydrodynamic limits~\cite{2016brownian}.

Afterwards, these authors joined by Sergio Simonella have also studied the strict linearization of the Boltzmann equation, also proving results on long time scales~\cite{2023longcor, 2024longder}.

The present paper is dedicated to an improvement of the convergence rate for the Rayleigh gas model~\cite{2016brownian}, which represents more than an exponential gain, along with the correction of a few inaccuracies, some of which having already been the subject of an erratum to the 2013 paper. Here, to obtain this better convergence rate, the main idea happens in the pruning process: the dynamics is written in terms of collision trees that are truncated at specific times so as to get reasonable bounds, and the error due to this truncation is the biggest one appearing in the proof. The original article used to cut these trees at regular intervals, but choosing a dividing more adapted to the situation greatly reduces this error. This kind of idea is frequent in the literature: compare for example the uniform cutting of the frequency space in the article of Hajer Bahouri and Jean-Yves Chemin~\cite{1999bahourichemin} and the more adapted cutting by Daniel Tataru~\cite[Figure~1]{2002tataru}, yielding keener Strichartz estimates for a quasilinear wave equation. Hence, for our problem, while the previous convergence result could hardly have a physical meaning -- as it yielded a rate in $(\log \log N)^{-1}$ -- the new one is much more satisfying with a rate approaching $\exp(-c \log N) = N^{-c}$.

Furthermore, while this rate was degrading with the time horizon in the original article~\cite{2016brownian}, it is now completely independent of the time parameter, as soon as this horizon~$t$ is chosen small enough compared with the number of particles~$N$. This dependence between $t$ and $N$ appears so as to handle the pruning process, and leads us to impose the same hydrodynamic scaling as in~\cite{2016brownian}: we only improve the kinetic scaling.

The main theorem is stated in Section~\ref{sec:res} and its proof is the subject of Section~\ref{sec:adaptpru}, first precising one of the original paper's estimates, then getting a better bound on the pruned-out term, and eventually concluding by tuning finely the pruning parameters' scaling. Section~\ref{sec:model} starts by recalling the framework and preliminary results.

\section{Modelling a tagged particle in a gas at equilibrium} \label{sec:model}

\subsection{Microscopic model}

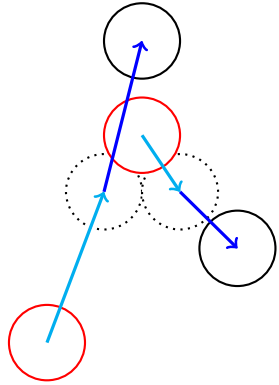
\begin{figure}[ht] 
\centering 
\begin{tikzpicture}
\newcommand{\scale}{1/2}
\draw[black, thick, dotted, fill = white] (0,0) circle (\scale*1 cm);
\draw[black, thick, dotted, fill = white] (2*\scale,0) circle (\scale*1 cm);
\draw[red, thick, fill = white] (\scale,\scale*1.5) circle (\scale*1 cm);
\draw[red, thick, fill = white] (-1.5*\scale,-4*\scale) circle (\scale*1 cm);
\draw[black, thick, fill = white] (\scale,4*\scale) circle (\scale*1 cm);
\draw[black, thick, fill = white] (3.5\scale,-\scale*1.5) circle (\scale*1 cm);
\draw[cyan, very thick, ->] (-1.5*\scale,-4*\scale) -- (0,0);
\draw[blue, very thick, ->] (2*\scale,0) -- (3.5\scale,-\scale*1.5);
\draw[cyan, very thick, ->] (\scale,\scale*1.5) -- (2*\scale,0);
\draw[blue, very thick, ->] (0,0) -- (\scale,4*\scale);
\end{tikzpicture}
\caption{Positions slightly before (in red) and after (in black) a collision}
\label{fig:coll1}
\end{figure}

The state of the gas of $N$ particles that we study is completely determined by the position (in the $d$-dimensional torus $\T^d$) and the velocity of each particle, represented by the vector
\begin{equation}
\uz_N = (z_1, \dots, z_N) \doteq (\ux_N, \uv_N) \in \mD^N \doteq (\T^{d} \times \R^{d})^N.
\end{equation}
The \emph{hard sphere} model consists in an exclusion condition, which states that two particles cannot get closer than a certain diameter~$\varepsilon$ -- so that we work with the following restricted open domain
\begin{equation} \label{eq:dom}
\mathcal{D}^\varepsilon_N = \{ \uz_N \in \mathcal{D}^N;\ \forall\ i\neq j,\ |x_i - x_j| > \varepsilon \}.
\end{equation}
Within $\mathcal{D}^\varepsilon_N$, the particles' dynamics is given by Newton equations for uniform line movement, i.e.
\begin{equation}
\frac{\dd x_i}{\dd t } = v_i, \hspace{15mm}
\frac{\dd v_i}{\dd t } = 0.
\end{equation}
Conversely, on the boundary of $\mathcal{D}^\varepsilon_N$, at least two particles are in contact. For instance, for the collision of the pair of particles $(i,j)$, we hence have $|x_i - x_j| = \varepsilon$. In this hard sphere model, the interaction is instantaneous and elastic. If the scalar product $(x_i - x_j) \cdot (v_i - v_j)$ is positive, the uniform line movement is well-defined, but otherwise the post-collisional velocities $({v_i}', {v_j}')$ are given by the following system, as pictured in Figure~\ref{fig:coll1},
\begin{equation} \label{eq:coll}
\left\{ \begin{array}{l}
{v_i}' = v_i - \frac{1}{\varepsilon^2} \left[ (v_i - v_j)\cdot(x_i - x_j)\right](x_i - x_j) \\
\ \\
{v_j}' = v_j + \frac{1}{\varepsilon^2} \left[(v_i - v_j)\cdot(x_i - x_j)\right](x_i - x_j).
\end{array} \right.
\end{equation}
Finally, we consider $f_N(t,~\cdot~)$ the probability density of presence of particles on the phase space $\mD^\varepsilon_N$ at time $t\geq 0$, which must be symmetrical by exchangeability of the particles. The microscopic dynamics provides the Liouville transport equation for $f_N$ on $\mD^\varepsilon_N$ 
\begin{equation} \label{eq:liouv}
\partial_t f_N + \uv_N\cdot \nabla_{\ux_N}  f_N = 0,
\end{equation}
with the following boundary condition on the post-collisional states:
\begin{equation}
|x_i - x_j|=\varepsilon \und (x_i - x_j)\cdot(v_i - v_j) > 0 \ \ \ \ \Rightarrow\ \ \ \  f_N(\uz_N) \doteq f_N(\uz_N^\star), \label{eq:boundcond}
\end{equation}
where $\uz_N^\star = (z_1, \dots, x_i, {v_i}^\star, \dots,  x_j, {v_j}^\star, \dots, z_N)$ denotes the \emph{pre}-collisional state associated to $\uz_N$. The dynamics is well-defined in these terms, up to a set of initial configurations of measure zero, as proved by another one of Lanford's students, Roger Keith Alexander~\cite{1975alexander, 2013newton}.

Other models implying different potentials of interaction have been studied, short-range~\cite{2013newton} or long-range~\cite{1999BoltLongrange, 2017linnoncutoff}. For a complete review of collisional kinetic theory, see the one by C\'{e}dric Villani~\cite{2002villanoche}.

\subsection{Nonideal Rayleigh gas and linear Boltzmann equation}

The kinetic limit we consider is called the low density limit, or Boltzmann--Grad limit, and consists in letting the number of particles $N$ go to infinity while keeping a constant mean free path $N^{-1} \varepsilon^{1-d} = 1$, so that the particles' diameter~$\varepsilon$ goes to 0.
In this limit, assuming initial chaos, the previously exposed framework usually leads to the Boltzmann equation, at least for short times~\cite{1975lanford, 2013newton}. 

In this paper, like in the original article~\cite{2016brownian}, we choose initial conditions close to equilibium to retrieve a linear version of the Boltzmann equation, whose theory is much simpler and hence might be derived for long time scales. More precisely, we consider the \emph{nonideal Rayleigh gas} model~\cite{1988spohnkineq}: we choose to tag the first particle, breaking the symmetry that was previously stemming from the particles' exchangeability, and we consider the initial condition
\begin{equation}
f_N(0, \uz_N) = \frac{\ind_{\mD^\varepsilon_N}(\uz_N)}{\mZ_N} \rho(x_1) M_\beta^{\otimes N}(\uz_N),
\end{equation}
where $\rho$ is a continuous space perturbation on the torus, $\mZ_N$ is a normalization constant, and $M_\beta$ denotes the equilibrium Maxwell state
\begin{equation} \label{def:max}
M_\beta(x,v) \doteq \left(\frac{\beta}{2\pi} \right)^{d/2} \exp\left( -\frac{\beta}{2} |v|^2 \right).
\end{equation}
The parameter~$\beta$ stands for an inverse temperature, tuning the intensive (kinetic) energy of the system. We take interest in the marginals of $f_N$, defined for $n\leq N$ as
\begin{equation}
f_N^{(n)}(t, \uz_n) \doteq \int_{\mD^{N-n}} f_N(t, \uz_n, z_{n+1}, \dots, z_N) \ind_{\mD_N^\varepsilon}(\uz_N) \dd z_{n+1} \dots \dd z_N. \label{eq:marg}
\end{equation}
Then, in the Boltzmann--Grad limit, the first marginal $f_N^{(1)}$ behaves like the solution $g\doteq M_\beta \varphi$ to the linear \emph{Rayleigh--Boltzmann equation}~\cite{2016brownian}, where
\begin{equation} \label{eq:phi}
\partial_t \varphi + v\cdot \nabla_x \varphi =  \int_{\Sf^{d-1}} \int_{\R^d} [\varphi(v^\star) - \varphi(v)] M_\beta(v_c) \left[ \omega \cdot (v_c - v) \right]_+ \dd v_c \dd \omega,
\end{equation}
with initial condition
\begin{equation} \label{eq:phiCI}
\varphi(0,x,v) = \rho(x). 
\end{equation}
Throughout this paper, we simply denote $\|\rho\| \doteq \| \rho \|_{\Lp^\infty(\T^d)}$.
This linear equation~(\ref{eq:phi}, \ref{eq:phiCI}) is globally well-posed in $\Lp^\infty(\T^d \times \R^d)$, and becomes a linear heat equation in the hydrodynamic limit~\cite{2016brownian}.

Some partial results exist for this same model with long-range interactions instead of hard sphere collisions~\cite{1999BoltLongrange, 2017linnoncutoff}, yet the complete derivation of the Rayleigh--Boltzmann~equation for general potentials is still an open problem. Other ways to derive the linear Rayleigh--Boltzmann~equation for long time scales are the \emph{ideal} Rayleigh gas model, in which the particles at equilibrium don't interact among themselves~\cite{1992centralray, 2018rayleigh, 2019rayleighannil}, and the Lorentz gas model, which consists in letting a tagged particle evolve in a frozen random background~\cite{1991spohnintpart, 1994CIP, 2008golselorentz}.

\section{Improved convergence rate of the first marginal} \label{sec:res}
As announced in the previous sections, the following theorem provides a convergence rate of the density's first marginal -- which corresponds to the tagged particle -- to the solution of the Rayleigh--Boltzmann equation. This is exactly the same convergence as in~\cite{2016brownian}, yet the rate depending on the number of particles $N = \varepsilon^{1-d}$ has been improved, from an error of order $(\log \log N)^{-1}$ to an error of order  $\exp\left( - c_\beta \left\lvert \log N \right\rvert^{1-\alpha} \right)$ for any $\alpha >0$, which is an improvement by more than an exponential.

\begin{theo}[Convergence to the Rayleigh--Boltzmann density] \label{theo1}
There exists a constant $c_\beta$ depending only on the temperature and the dimension such that, for any $\alpha \in (0,1/2)$, as long as 
\begin{equation} t \lesssim \left( \log\left| c_\beta \log \varepsilon \right|\right) ^{\frac{1}{2} - \alpha}, 
\end{equation}
one has -- for $\varepsilon$ small enough -- the following convergence rate of the BBGKY distribution's first marginal to the linear Boltzmann density, in the low density limit $N = \varepsilon^{-(d-1)} \rightarrow \infty$,
\begin{equation}
\lnorm f_N^{(1)} - g \rnorm_{\Lp^\infty([0,t]\times \mD^d)} \leq \| \rho \| \exp\left( - c_\beta \left\lvert \log \varepsilon \right\rvert^{1-\alpha} \right).
\end{equation}
\end{theo}
The notation $t \lesssim \left( \log\left| c_\beta \log \varepsilon \right|\right) ^{\frac{1}{2} - \alpha}$ means that for a constant $c$ small enough one has 
\begin{equation}
t \leq c \left( \log\left| c_\beta \log \varepsilon \right|\right) ^{\frac{1}{2} - \alpha}.
\end{equation}

The proof of this theorem is the subject of Section~\ref{sec:adaptpru}. The main idea happens in the pruning process, which consists in removing the trajectories containing too many collisions at certain intermediate times. Note that although we managed to get rid of the time dependence in the convergence rate, the time scaling does not get better than in the original article~\cite{2016brownian} and remains of order $\left( \log\left| c_\beta \log \varepsilon \right|\right) ^{\frac{1}{2} - \alpha}, $ so that the scaling of the hydrodynamic limit in~\cite{2016brownian} does not get improved either. Indeed, one may see in the pruning process (Proposition~\ref{prop:estrem}) that the length of the time horizon~$t$ is deeply coupled to the pruning parameter~$K$ which cannot get too large compared to~$N$ (see~Proposition~\ref{prop:estpru}).

\section{Adaptative pruning and proof of the convergence rate} \label{sec:adaptpru}

This section is dedicated to the presentation of a tree pruning tuned adaptatively in time, so as to improve the convergence rate which could be found in the original article~\cite{2016brownian}, as stated in Theorem~\ref{theo1}. 
Let us recall the notation and general framework of our study. We have the following hierarchy on the marginals of the density~$f_N$, called \emph{BBGKY hierarchy} and defined in Boltzmann's lectures~\cite{1872boltzmann} , in the form of a Dyson expansion
\begin{equation} \label{eq:Dys}
f_N^{(n)}(t) = \sum_{s=0}^{N-n} Q_{n,n+s}(t) f_N^{(n+s)}(0),
\end{equation}
where the \emph{successive-collision operators} are given by
\begin{equation} \label{def:Q}
Q_{n,n+s}(t) \doteq  \int_0^t \int_0^{t_1} \cdots \int_0^{t_{s-1}} \Theta_n(t-t_1)\mathcal{C}_n \Theta_{n+1}(t_1 - t_2)\mathcal{C}_{n+1}\dots \Theta_{n+s}(t_s) \dd t_s \dots \dd t_{1},
\end{equation}
and are written in terms of $\Theta_n$, the operator associated to the free transport with collisions in $\mD^\varepsilon_n$ with specular reflection (see for example~\cite{2002villanoche} for details on this operator's definition), and in terms of the elementary \emph{collision operators} defined as follows
\begin{equation} \label{def:Cgrandcan}
C_n f_N^{(n+1)}(\uz_n) \doteq (N-n)\varepsilon^{d-1} \sum_{i=1}^n \int_{\R^d}\int_{\Sf^{d-1}} \omega\cdot (v_{n+1} - v_i) f_N^{(n+1)}(\uz_n, x_i + \varepsilon \omega, v_{n+1}) \dd \omega \dd v_{n+1}.
\end{equation}
All of these operators have a formal limit version in the Boltzmann--Grad limit, where $(N-n)\varepsilon^{d-1}$ converges to 1, while $\varepsilon$ goes to 0. They easily satisfy the same estimates as their BBGKY counterparts, estimates that are presented in the following sections.

\subsection{A continuity estimate on the successive-collision operators}
We introduce the total kinetic energy
\begin{equation}
H_{k}( \uv_k) = \frac{1}{2}\sum_{i = 1}^k |v_i|^2,
\end{equation}
and for $\lambda > 0$ holding the role of an inverse temperature and $k\in \N^*$, we consider the space $\mF_{\varepsilon, k, \lambda}$ of measurable functions defined almost everywhere on the restricted domain $\mD^\varepsilon_k$ such that
\begin{equation} \label{def:foncspace}
\|f_k\|_{\varepsilon, k, \lambda} \doteq \supess_{\uz_k \in \mD^\varepsilon_k} \Bigl |f_k(\uz_k) \exp(\lambda H_k(\uv_k)) \Bigr | < \infty.
\end{equation}
The marginals of $f_N$ belong to this space for $\lambda = \beta$, as stated in~\cite{2016brownian}, with
\begin{equation} \label{ineq:Fwnorm}
\sup_{t \geqslant 0} \ \| f_N^{(k)}(t) \|_{\varepsilon, k, \beta } \leq \sup_{\uz_k \in \mD^\varepsilon_k} M_\beta^{\otimes k}(\uv_k) \exp\bigl(\beta H_k(\uv_k)\bigr) \|\rho\| = \|\rho\| \left(\frac{\beta}{2\pi}\right)^{\frac{kd}{2}}\cdotp
\end{equation}
This bound is very strong and specific to the linear case: it is valid for all times and this is the reason why we can derive linear results on long time scales; in the non-linear case the a~priori bounds were only valid for very small times~\cite{2013newton}.

The following proposition is very similar to its equivalent in the original article~\cite[Lemma~4.2]{2016brownian}; it is simply a little bit more general concerning the degrading rate of the norm: indeed, the idea is to show a continuity estimate between two of the spaces defined above for different temperature parameters, degrading this parameter so as to resorb a factor~$|\uv|$ by a fraction of the sub-Gaussian decreasing at infinity, which is precisely tuned by the temperature in the considered spaces.
This continuity estimate, along with the bound on the marginals~(\ref{ineq:Fwnorm}), justify the convergence of the Dyson expansion~(\ref{eq:Dys}) in $\mF_{\varepsilon, s, \beta/2}$ for times $t \geq 0$ small enough.

\begin{prop}[Continuity of the successive-collision operators] \label{prop:contestQ}\ \\
There exists a constant $C_d$ depending only on the dimension such that for all~$n,s \in \N^*$ and all times~$t>0$, if $f_{n+s} \in \mF_{\varepsilon,n+s,\lambda}$, then for every $b \geq 2$, we have
\begin{equation}
\begin{array}{c}
Q_{n, n+s}(t) f_{n+s} \in \mF_{\varepsilon, n, \lambda \left(1-\frac{1}{b}\right)} \textit{, with} \\ \ \\
\Bigl \| Q_{n, n+s}(t) f_{n+s} \Bigr \|_{\varepsilon, n, \lambda \left(1-\frac{1}{b}\right)} \leq e^{n} \left( \frac{ \sqrt{b} \ C_d t}{\lambda^{(d+1)/2}} \right)^s \|f_{n+s}\|_{\varepsilon, n+s, \lambda}.
\end{array}
\end{equation} 
\end{prop}

\begin{prof}
First of all, let us observe that the transport operators preserve all of the weighted norms $\left(\|~\cdot~\|_{\varepsilon, n, \lambda}\right)_{n,\lambda}$, since the weight depends only on the kinetic energy of the system.

Then, for $j\leq N$,  using that $(N-j)\varepsilon^{d-1} \leq 1$ in our scaling, let us compute for $f_{j+1} \in \mF_{\varepsilon, j+1, \lambda}$, making its norm appear,
\begin{align*}
\Bigl| \Theta_{j}(-\tau) \mC_{j} \Theta_{j+1}(\tau) f_{j+1}   \Bigr| & \leq \left| \Theta_{j}(-\tau) \sum_{i=1}^j \int \omega\cdot (v_{j+1} - v_i) \Theta_{j+1}(\tau) f_{j+1}(\uz_j, x_i + \varepsilon \omega, v_{j+1}) \dd \omega \dd v_{j+1} \right| \\
& \leq \sum_{i=1}^j \int_{\Sf^{d-1} \times \R^d} (|v_{j+1}| + |v_i|) \|f_{j+1}\|_{\varepsilon,j+1, \lambda} \exp\left[-\lambda H_{j+1}(\uv_{j+1}) \right] \dd \omega \dd v_{j+1}\\
& \ \ = |\Sf^{d-1}|\cdot \|f_{j+1}\|_{\varepsilon,j+1, \lambda} \sum_{i=1}^j \int_{\R^d} (|v_{j+1}| + |v_i|)  \exp\left[-\frac{\lambda}{2} \sum_{k=1}^{j+1} |v_k|^2 \right] \dd v_{j+1}.
\end{align*}
The latter integrals may be written explicitly, up to constants depending only on the dimension~$d$, after a radial change of variable and a dilation by $\lambda^{-\frac{1}{2}}$,
\begin{align}
\int (|v_{j+1}| + |v_i|)  \exp\left[-\frac{\lambda}{2} \sum_{k=1}^{j+1} |v_k|^2) \right] \dd v_{j+1} & = C_d \int (r + |v_i|) \exp\left[-\frac{\lambda}{2} \sum_{k=1}^{j} |v_k|^2 \right] r^{d-1} e^{-\frac{\lambda}{2} r^2} \dd r \nonumber \\
& = C_d \exp\Bigl[-\lambda H_j(\uv_j)\Bigr] \left( c_d \sqrt{\lambda^{-(d+1)}} + |v_i| \tilde{c}_d \sqrt{\lambda^{-d}}\right).	 \label{eq134}
\end{align}
This way, applying this to $f_{n+s}$, summing~(\ref{eq134}) over $i$ and accepting to downgrade by $\lambda / bs$ the considered norm so as to later resorb the factors $|v_i|$, we get that 
\begin{align} \label{eq23}
\begin{array}{l}
\Bigl \lVert \Theta_{n+s-1}(-t_s)\mathcal{C}_{n+s-1} \Theta_{n+s}(t_s)f_{n+s}\Bigr \rVert_{\varepsilon, n+s-1, \lambda  - \lambda / bs} \\
 \hspace{1cm}\leq \widetilde{C}_d \left( (n+s-1) \sqrt{\lambda^{-(d+1)}} + \sqrt{\lambda^{-d}} \sum_{i=1}^{n+s-1} |v_i| \right) \exp\left[-\frac{\lambda}{2bs} \sum_{k=1}^{n+s-1} |v_k|^2 \right] \|f_{n+s}\|_{\varepsilon,n+s, \lambda}. 
\end{array}
\end{align}
But using the Cauchy--Schwarz inequality and the fact that $x e^{-x} \leq e^{-1} $ for any $x \geq 0$ we have
\begin{align}
\left( \sum_{i=1}^{n+s-1} |v_i| \right) \exp\left[-\frac{\lambda}{2bs} \sum_{k=1}^{n+s-1} |v_k|^2 \right] & \leq \left(  \frac{(n+s-1) bs}{\lambda} \right)^{\frac{1}{2}} \left( \sum_{i=1}^{n+s-1} |v_i|^2 \frac{\lambda}{bs} \right)^{\frac{1}{2}}  e^{-\frac{\lambda}{2bs} \sum_{k=1}^{n+s-1} |v_k|^2} \nonumber \\
&  \leq  (n+s) \sqrt{\frac{b}{ \lambda e}}\ , \label{eq:CS}
\end{align}
so that~(\ref{eq23}) yields
\begin{equation} \label{eq:theta}
\Bigl \lVert \Theta_{n+s-1}(-t_s)\mathcal{C}_{n+s-1} \Theta_{n+s}(t_s)f_{n+s}\Bigr \rVert_{\varepsilon, n+s-1, \lambda  - \lambda / bs} \leq \widehat{C}_d (n+s) \sqrt{b}  \lambda^{-\frac{d+1}{2}} \|f_{n+s}\|_{\varepsilon,n+s, \lambda}. 
\end{equation}
To retrieve $Q_{n, n+s}(t)$, we have to iterate this calculus $s$ times, downgrading the norm parameter by a factor $(1 - 1 / bs)$ at each step. At the $i$-th iteration, the factor in the right-hand side of~(\ref{eq:theta}) becomes 
\begin{align}
\widehat{C}_d (n+s) \sqrt{b}  [\lambda(1-1/bs)^{i-1}]^{-\frac{d+1}{2}} & \leq \widehat{C}_d (n+s) \sqrt{b}  [\lambda(1-1/bs)^s]^{-\frac{d+1}{2}}. 
\end{align}
Using the convexity of $x \mapsto (1 - x)^s$, and then the fact that $b \geq 2$, we get that
\begin{align}
\widehat{C}_d (n+s) \sqrt{b}  [\lambda(1-1/bs)^{i-1}]^{-\frac{d+1}{2}} & \leq \widehat{C}_d (n+s) \sqrt{b}  [\lambda(1-1/b)]^{-\frac{d+1}{2}} \label{lignconv}\\
& \leq \widehat{C}_d (n+s) \sqrt{b}  [\lambda / 2]^{-\frac{d+1}{2}}. \label{lignb}
\end{align}
The final inverse temperature parameter, which is smaller than all the intermediate ones, is
\begin{equation}
\lambda\left(1 - \frac{1}{bs}\right)^{s} \geq \lambda\left(1 - \frac{1}{b}\right),
\end{equation}
by the same convexity argument, so that eventually, since the considered norms are increasing with~$\lambda$,
\begin{align*}
\Bigl \lVert Q_{n, n+s}(t) f_{n+s} \Bigr \rVert_{\varepsilon, n, \lambda - \frac{\lambda }{ b}}& \leq \left(\overline{C}_d (n+s) \sqrt{b \lambda^{-(d+1)}} \right)^s \int_{0 \leq t_s \leq \dots \leq t_1 \leq t}  \dd t_1 \cdots \dd t_s \|f_{n+s} \|_{\varepsilon, n+s, \lambda}\\
&\leq \left(\overline{C}_d  \sqrt{b \lambda^{-(d+1)}} \right)^s \left(n+s\right)^s  \ \frac{t^s}{s!}\ \|f_{n+s} \|_{\varepsilon, n+s, \lambda} \\
& \leq \left(\overline{C}_d  \sqrt{b \lambda^{-(d+1)}} \right)^s e^{n+s} \ t^s\ \|f_{n+s} \|_{\varepsilon, n+s, \lambda}\ ,
\end{align*}
where the factor $s!$ comes from the imposed order of collision times~$t_s \leq \dots \leq  t_1$, and allows to control the term $(n+s)^s$, concluding the proof. \CQFD
\end{prof}

\subsection{Pseudo-trajectories' adaptative tree pruning}

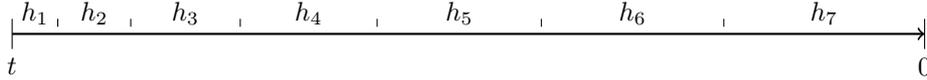
\begin{figure}[h!] 
\centering
\begin{tikzpicture}
\newcommand{\ho}{0.05}
\newcommand{\hd}{0.13}
\newcommand{\htr}{0.25}
\newcommand{\hc}{0.4}
\newcommand{\hp}{0.58}
\newcommand{\hs}{0.78}
\draw[black, thick, ->] (0,0) -- (12,0);
\draw[black] (12,0.2) -- (12,-0.2) node[below]{$0$};
\draw[black, thin, dashed] (\ho*12,0.2) -- (\ho*12,0);
\draw[black, thin, dashed] (\hd*12,0.2) -- (\hd*12,0) ;
\draw[black, thin, dashed] (\htr*12,0.2) -- (\htr*12,0);
\draw[black, thin, dashed] (\hc*12,0.2) -- (\hc*12,0);
\draw[black, thin, dashed] (\hp*12,0.2) -- (\hp*12,0);
\draw[black, thin, dashed] (\hs*12,0.2) -- (\hs*12,0);
\draw[black, thin, dashed] (\ho*6,0) node[above]{$h_1$};
\draw[black, thin, dashed] (\hd*6 + \ho*6,0) node[above]{$h_2$};
\draw[black, thin, dashed] (\htr*6+\hd*6,0) node[above]{$h_3$};
\draw[black, thin, dashed] (\hc*6 + \htr*6,0) node[above]{$h_4$};
\draw[black, thin, dashed] (\hp*6 + \hc*6,0) node[above]{$h_5$};
\draw[black, thin, dashed] (\hs*6 + \hp*6 ,0) node[above]{$h_6$};
\draw[black, thin, dashed] (\hs*6 + 6 ,0) node[above]{$h_7$};
\draw[black] (0, 0.2) -- (0,-0.2) node[below]{$t$};
\end{tikzpicture}
\caption{Backwards division of the time interval under study}
\label{fig4}
\end{figure}

Now, so as to work with collision trees of controlled sizes, we will simply consider truncated series instead of our Dyson expansions~(\ref{eq:Dysexprem}). More precisely, for a fixed time $t>1$ we will study how our functionals behave on a well chosen cutting of the interval $[0,t]$, and impose a maximal amount of collisions on each small piece of this cutting. We thus divide the interval $[0,t]$ into $K \in \N^*$ little pieces of sizes $h_1, \dots, h_K$ respectively, backwards in time, starting from time $t$ and going back to zero, as pictured in Figure~\ref{fig4}.

We will morally forbid more than 2 collisions per particle to happen in each small interval: at the $k$-th time quantum of length $h_k$, we want that at most $2^k$ particles may have collided, as if we were pruning the collision tree every time it becomes more than exponentially big. To get an explicit formulation of this condition on the series expansion, we will write the Dyson series~(\ref{eq:Dys}) between $t$ and $t-h_1$, cut it after 2 collisions, then do it again between $t-h_1$ and $t-h_2$ after $2^2$ collisions, and eventually iterate this calculus $K$ times: we denote the time steps
\begin{equation}
t^\mathrm{p}_k = t - \sum_{j=1}^k h_i,
\end{equation}
with the condition \begin{equation}
\sum_{i=1}^K h_i = t,
\end{equation}
so that $t_K^\mathrm{p} = 0$. So we write as in~(\ref{eq:Dys})
\begin{align}
& f_N^{(1)}(t) =  \sum_{j_1 = 0}^{1}  Q_{1, 1+j_1}(h_1)f_N^{(1 + j_1)}(t-h_1) + \sum_{s = 2}^\infty Q_{1, 1+s}(h_1)f_N^{(1 + s)}(t-h_1) \\
&= \sum_{j_1 = 0}^{1} \cdots \sum_{j_K = 0}^{2^K - 1} Q_{1, J_1}(h_1) Q_{J_1, J_2}(h_2) \dots Q_{J_{K-1}, J_K}(h_K) f_N^{(J_K)}(0) \label{eq:Dysexprem}\\
& +  \sum_{k=1}^K \sum_{j_1 = 0}^{1} \cdots \sum_{j_{k-1} = 0}^{2^{k-1} - 1}  Q_{1, J_1}(h_1) \dots  Q_{J_{k-2}, J_{k-1}}(h_{k-1})  \sum_{s=2^k}^\infty Q_{J_{k-1}, J_{k-1} + s}(h_{k}) f_N^{(J_{k-1} + s)}(t^{\mathrm{p}}_k),\label{eq:Dysexpru} 
\end{align}
denoting
\begin{equation}
J_K \doteq 1 + j_1 + \cdots + j_K.
\end{equation}
We hence introduce the following remainder, which corresponds to the pruned-out trajectories,
\begin{equation} \label{eq:remK}
R^{[K]}(t) \doteq \sum_{k=1}^K \sum_{j_1 = 0}^{1} \cdots \sum_{j_{k-1} = 0}^{2^{k-1} - 1}  Q_{1, J_1}(h_1) \dots  Q_{J_{k-2}, J_{k-1}}(h_{k-1})  \sum_{s=2^k}^\infty Q_{J_{k-1}, J_{k-1} + s}(h_{k}) f_N^{(J_{k-1} + s)}(t^{\mathrm{p}}_k) .
\end{equation}
We will give estimates on the truncated series in the following section, but first we have to justify that the remainder is small enough. This is the point of the following proposition, using the continuity estimates of previous section, and improving greatly the results of~\cite[Proposition~4.3]{2016brownian} by adapting the time cutting $(h_1, \dots, h_K)$. Since the chosen condition of a sub-exponential number of collisions is very restrictive at first, and then gradually relaxed, the key point is to chose the cutting times small at first and then progressively bigger.

\begin{prop}[Estimate of the pruned-out term] \label{prop:estrem}
With the previous notation, for any $\alpha \in (0,1/2)$ and $K$ large enough satisfying $t \lesssim K^{\frac{1}{2}-\alpha}$, a good choice of time cutting $\uh = (h_1, \dots, h_K)$ provides the following estimate 
\begin{equation}
\lnorm R^{[K]}(t) \rnorm_{\Lp^\infty(\mD^d)}  \leq \| \rho \| e^{-2^{K - K^\alpha}}. 
\end{equation}
\end{prop}

\begin{prof} As $f_N^{(J_{k-1} + s)} \in\mF_{\varepsilon, J_{k-1} + s, \beta}$ by the bound~(\ref{ineq:Fwnorm}) on its norm, the continuity estimate on the successive-collision operators  given in Proposition~\ref{prop:contestQ} asserts that for any $h_k$ small enough, 
\[ \begin{array}{l}
\sum_{s=2^k}^\infty Q_{J_{k-1}, J_{k-1} + s}(h_k) f_N^{(J_{k-1} + s)}(t^\mathrm{p}_k)  \in \mF_{\varepsilon, J_{k-1}, \beta /2}
\end{array} \]
with
\begin{align*} \def\arraystretch{3} 
\begin{array}{l}
\lnorm \sum_{s=2^k}^\infty Q_{J_{k-1}, J_{k-1} + s}(h_k) f_N^{(J_{k-1} + s)}(t^\mathrm{p}_k) \rnorm_{\varepsilon, J_{k-1}, \beta /2}\\
\hspace*{30mm} \leq e^{J_{k-1}} \sum_{s=2^k}^\infty \left( \frac{ \sqrt{2} C_d  h_k}{\beta^{(d+1)/2}} \right)^s \lnorm  f_N^{(J_{k-1} + s)}(t^\mathrm{p}_k)\rnorm_{\varepsilon, J_{k-1}+s, \beta}\\
\hspace*{30mm} \leq e^{J_{k-1}} \sum_{s=2^k}^\infty \left( \frac{ \sqrt{2} C_d  h_k}{\beta^{(d+1)/2}} \right)^s (C \beta^{d/2})^{J_{k-1}+s} \| \rho \|.
\end{array}
\end{align*}
We now iterate $k$~times Proposition~\ref{prop:contestQ} -- like in the proof of this same proposition, downgrading each time the norm by a factor $1-1/2k$ so that the final inverse temperature, which is smaller than all the intermediate ones, can be bounded in the following way by convexity of $x \mapsto (1-x)^k$:
\begin{equation}
\frac{\beta}{2} \left( 1 - \frac{1}{2k} \right)^k \geq \frac{\beta}{2} \left( 1 - \frac{1}{2} \right) = \frac{\beta}{4}\cdotp
\end{equation}
This way, our $k$ iterations of Proposition~\ref{prop:contestQ} -- choosing $b=2k$ -- allow us to write, grouping all the terms appearing under the form $e^{J_i}$ or $(C\beta^{d/2})^{J_{k-1}}$ together as a power of a constant $C(\beta)$,
\[ \begin{array}{l}
\lnorm Q_{1, J_1}(h_1) \dots  Q_{J_{k-2}, J_{k-1}}(h_{k-1}) \sum_{s=2^k}^\infty Q_{J_{k-1}, J_{k-1} + s}(h_k) F_{J_{k-1} + s}(t^\mathrm{p}_k) \rnorm_{\varepsilon, 1, \beta/4} \\  
\hspace*{3mm} \leq \| \rho\|  C(\beta)^{\sum_{i=0}^{k-1} J_i} \left( \sqrt{2k} C_d  h_1 \frac{4^{(d+1)/2}}{\beta^{(d+1)/2}} \right)^{j_1} \dots  \left( \sqrt{2k} C_d  h_{k-1} \frac{4^{(d+1)/2}}{\beta^{(d+1)/2}} \right)^{j_{k-1}} \sum_{s=2^k}^{\infty} \left( \frac{\sqrt{2} C_d  h_k}{\sqrt{\beta}} \right)^s \cdotp
\end{array} \]
Let us then observe that since $J_i = 1 + j_1 + \dots + j_i$, we have the following bound on
\[ \sum_{i=0}^{k-1} J_i \leq \sum_{i=0}^{k-1} (1 + 2 + \dots + 2^i) \leq \sum_{i=0}^{k-1} 2^{i+1} \leq 2^{k+1}. \] 
Hence, since all of the weighted norms are greater than the $\Lp^\infty$-norm, up to a new constant~$C$ depending on $\beta$ we can write
\begin{equation}\label{eq:estRR} \begin{array}{l}
\lnorm R^{[K]}(t) \rnorm_{\Lp^\infty}  \leq \|\rho\|   \sum_{k=1}^K C^{2^{k+1}} \sum_{j_1 = 0}^{1} \cdots \sum_{j_{k-1} = 0}^{2^{k-1} - 1} \left( \sqrt{k} C h_1 \right)^{j_1} \dots  \left( \sqrt{k}  C  h_{k-1} \right)^{j_{k-1}} \sum_{s=2^k}^{\infty} \left(  C  h_k \right)^s. \\
\end{array} 
\end{equation}

Let us henceforth tune our cutting times $(h_1, \dots, h_K)$. We can see in the last equation that the last sum -- which will provide smallness -- will be of order $(C h_k)^{2^k}$, so that if the first cutting times need to be small, the following ones may get progressively bigger, as pictured in the following Figure~\ref{fig5}. Let us recall that this corresponds to the condition on a sub-exponential number of collisions being less and less restrictive.
\begin{figure}[h!] 
\centering
\scalebox{0.88}{
\begin{tikzpicture}
\newcommand{\longg}{16/42}
\newcommand{\ho}{1.2}
\newcommand{\hd}{3}
\newcommand{\htr}{6}
\newcommand{\hc}{11}
\newcommand{\hp}{17}
\newcommand{\hv}{24}
\draw[black, thick, ->] (0,0) -- (16,0);
\draw[black] (16,0.2) -- (16,-0.2) node[below]{$0$};
\draw[black, thin] (\ho*\longg,0.2) -- (\ho*\longg,0) ;
\draw[black, thin] (\ho*\longg/2,0.2)  node[above]{$h_1$} ;
\draw[black, thin] (\hv *\longg/2 +42*\longg/2,0.2)  node[above]{$h_K$} ;
\draw[black, thin] (\hp *\longg/2 + \hc *\longg/2 ,0.2)  node[above]{$\cdots$} ;
\draw[black, thin] (\hd*\longg,0.2) -- (\hd*\longg,0);
\draw[black, thin] (\htr*\longg,0.2) -- (\htr*\longg,0)  ;
\draw[black, thin] (\hc*\longg,0.2) -- (\hc*\longg,0) ;
\draw[black, thin] (\hv*\longg,0.2) -- (\hv*\longg,0) ;
\draw[black, thin] (\hd*\longg,0.2) -- (\hd*\longg,0) ;
\draw[black] (0, 0.2) -- (0,-0.2) node[below]{$t$};
\end{tikzpicture}}
\caption{Backwards construction of the cutting times}
\label{fig5}
\end{figure}
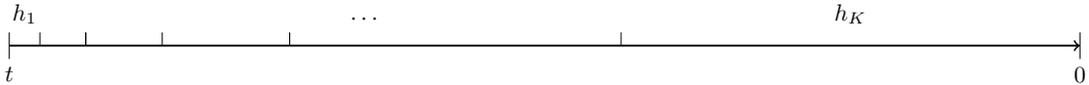

More precisely, we define for all $1 \leq i \leq K$,
\begin{equation}
\tilde{h}_i \doteq  \frac{e^{-2^{(K-K^{1- \alpha } - i)}}}{2 C \sqrt{K}},
\end{equation}
so that by positivity, and then the fact that $K - j  \geq K - K^{1-\alpha}$, hence $2^{K - K^{1-\alpha} - (K-j)} \leq 1,$  
\begin{align*}
\sum_{i = 1}^K \tilde{h}_i & \geq \sum_{j = 0}^{\lfloor K^{1-\alpha} \rfloor} \tilde{h}_{K-j} \\
& \geq \sum_{j = 0}^{\lfloor K^{1-\alpha} \rfloor}\frac{e^{-1}}{2 C \sqrt{K}} \\
& \geq t,
\end{align*}
since we required $t \lesssim K^{\frac{1}{2}-\alpha}$. This way, we might rescale this cutting into
\begin{equation}
h_i \doteq \frac{t}{\sum_{i = 1}^K \tilde{h}_i}\ \tilde{h}_i \leq \tilde{h}_i,
\end{equation}
so that \begin{equation} \sum_{i = 1}^K h_i = t. 
\end{equation}
Let us observe that the loss due to the $K^{1-\alpha}$ in the convergence rate is caused by the fact that we have to cover an interval of length $t$ that may go to infinity. For any finite value of $t$, this correction is not needed anymore and we may hence get a convergence rate in $e^{-2^K}$, which would give a better final estimate (see the final remark at page~\pageref{remark}).  \label{rem:finitime}

Eventually, for $K$ large enough, the estimate~(\ref{eq:estRR}) in the case of our cutting provides
\begin{align}
\lnorm R^{[K]}(t) \rnorm_{\Lp^\infty} &  \leq \|\rho\|  \sum_{k=1}^K C^{2^{k+1}} \prod_{i=1}^k \sum_{j_{i} = 0}^{2^{i} - 1} \left( \sqrt{k} C h_i \right)^{j_i}  \sum_{s=2^k}^{\infty} \left(  C  h_k \right)^s \\
&\leq \|\rho\| \sum_{k=1}^K C^{2^{k+1}} \prod_{i=1}^k \sum_{j_{i} = 0}^{2^{i} - 1} \left( \frac{1}{2} \right)^{j_i}  \sum_{s=2^k}^{\infty} \left( \frac{ e^{-2^{K-K^{1-\alpha}-k}}}{2\sqrt{K}} \right)^s,
\end{align}
so that bounding all of the $k$ geometric series by 2, then harnessing the factor $(\sqrt{K})^{-2^k}$ to crush the terms that blow up, we get as wanted
\begin{align}
\lnorm R^{[K]}(t) \rnorm_{\Lp^\infty} & \leq \| \rho \| \sum_{k=1}^K \left(C^{2} \right)^{2^k} \left[ \prod_{i=1}^k   2 \right]\times   2 \left( \frac{e^{-2^{K-K^{1-\alpha}-k}}}{2\sqrt{K}} \right)^{2^k}   \\
&\leq \| \rho \| \sum_{k=1}^K \left( \frac{2 C^{2}}{\sqrt{K}} \right)^{2^k} \left( e^{-2^{K-K^{1-\alpha}-k}}\right)^{2^k} \label{line:19}  \\
&\leq \| \rho \|    e^{-2^{K-K^{1-\alpha}}}, \nonumber
\end{align}
using at line~(\ref{line:19}) that $2^{k+1} \leq 2^{2^k}$, and hence concluding the proof.
\CQFD
\end{prof}

\subsection{Proof of the new convergence rate}
Let us recall that $g$ denotes the solution to the linear Rayleigh--Boltzmann equation given by~(\ref{eq:phi}) and~(\ref{eq:phiCI}). We denote  $R^{[K]}_{\lim}$ the pruned-out remainder of its Dyson expansion, exactly as for $f_N^{(1)}$ in~(\ref{eq:Dysexpru}), which satisfies easily the same estimate as the one given in Proposition~\ref{prop:estrem}. Now, the proximity of the pruned trajectories may be written in the following proposition.
\begin{prop}[Proximity of the pruned trajectories] \label{prop:estpru} There exists a constant $C_\beta$ depending on the temperature and on the dimension such that if $t \ll \sqrt{K}$, and if $K$ and $N = \varepsilon^{1-d}$ are large enough, then
\begin{equation}
\lnorm (f_N^{(1)} - R^{[K]}) - (g - R^{[K]}_{\lim} )\rnorm_{\Lp^\infty([0,t]\times \mD^d)} \leq {C_\beta}^{2^K} \|\rho\| \cdot  |\log \varepsilon|^{\frac{3d-1}{4}}\varepsilon^\frac{d-1}{2(d+1)}.
\end{equation}
\end{prop}

\begin{prof}
The complete proof of this proposition may be found in~\cite[Section~5]{2016brownian}, and follows from several approximations: an energy truncation and a time separation are operated so as to be able to construct a small set of bad collision parameters, which is such that out of this set one may consider that there is no recollision. Hence, the proof is brought back to studying the initial error at time $t=0$ and the very small error due to the prefactors~$(N-n)\varepsilon^{d-1}$.

The only improvement we bring to the original paper is the more precise estimate of the operators $Q_{1, J_1}(h_1) \dots  Q_{J_{k-1}, J_{k}}(h_{k})$, which like in the proof of Proposition~\ref{prop:estrem} induces a factor ${C_\beta}^{2^K}$, instead of the original crude bound with $|Q|_{1, J_{k}}(t)$, which gave a factor $(C_\beta t)^{2^K}$. Physically, we decompose the time interval into small  pieces whose lengths are adapted to the maximum number of particles that may appear in them, so that the dynamics behaves similarly during each one of them. Hence, as long as time does not get too big with respect to the number of pieces, none of the estimates depends on the total time length.

Also note that we have taken into account a correction in the geometric estimate given in~\cite[Proposition~5.1]{2016brownian}, which had been the subject of an erratum of~\cite{2013newton} and merely changes the power of~$\varepsilon$.

\CQFD
\end{prof}

We eventually obtain the following estimate of the convergence rate by tuning well our parameter $K$, this way proving  Theorem~\ref{theo1}.
\begin{prop}[Final convergence estimate] There exists a constant $c_\beta$ depending only on the temperature and the dimension such that, for any $\alpha \in (0,1/2)$, in the following scaling:
\begin{equation}
K = \left\lfloor \frac{\log\left( 2 c_\beta |\log \varepsilon | \right)}{\log 2}  \right\rfloor \und t \lesssim K^{\frac{1}{2} - \alpha},
\end{equation}
one has this final convergence rate of the BBGKY distribution to the Rayleigh--Boltzmann distribution
\begin{equation}
\lnorm f_N^{(1)} - g \rnorm_{\Lp^\infty([0,t]\times \mD^d)} \leq \| \rho \| \exp\left( - c_\beta \left\lvert \log \varepsilon \right\rvert^{1-\alpha} \right).
\end{equation}
\end{prop}

\begin{prof} Considering $C_\beta$ the constant given by Proposition~\ref{prop:estpru}, we choose 
\begin{equation}
K = \left\lfloor \frac{1}{\log 2} \log \left( \frac{ (d-1) |\log \varepsilon | }{4(d+1) \log C_\beta} \right) \right\rfloor,
\end{equation} 
so that we can choose the small constant as
\begin{equation}
c_\beta \doteq \frac{ (d-1) }{8(d+1) \log C_\beta}\cdotp
\end{equation}
In this scaling, one has precisely 
\begin{equation}
 {C_\beta}^{2^K} \leq \varepsilon^{\frac{1-d}{4(d+1)}},
\end{equation}
so that by~Proposition~\ref{prop:estpru},
\begin{equation}
\lnorm (f_N^{(1)} - R^{[K]}) - (g - R^{[K]}_{\lim} )\rnorm_{\Lp^\infty([0,t]\times \mD^d)} \leq \|\rho\|\cdot  |\log \varepsilon|^{\frac{3d-1}{4}} \varepsilon^{\frac{d-1}{4(d+1)}}.
\end{equation}
Hence, the only remaining error term is the remainder given by the pruning process, and for our scaling of $K$, Proposition~\ref{prop:estrem} yields that, since $K - K^\alpha \geq (1-\alpha) K$ for $K$ large enough,
\begin{align}
\lnorm R^{[K]}\rnorm_{\Lp^\infty([0,t]\times \mD^d)}&  \leq \| \rho \| e^{-2^{K - K^\alpha}} \\
& \leq \| \rho \| e^{-2^{(1-\alpha)K}},
\end{align} 
and eventually, writing 
\begin{equation}
K \geq \frac{\log\left( 2 c_\beta |\log \varepsilon | \right)}{\log 2} - 1 = \frac{1}{\log2} \log\left( c_\beta |\log \varepsilon | \right) ,
\end{equation}
we get
\begin{align}
\lnorm R^{[K]}\rnorm_{\Lp^\infty([0,t]\times \mD^d)}& \leq \| \rho \| \exp\left( - {c_\beta}^{1-\alpha} \left\lvert \log \varepsilon \right\rvert^{1-\alpha} \right)\\
& \leq \| \rho \| \exp\left( - \min(c_\beta, 1) \left\lvert \log \varepsilon \right\rvert^{1-\alpha} \right), 
\end{align} 
which is the biggest of both errors, concluding the proof.
\CQFD
\end{prof}
\hspace*{10mm}
\begin{rema} \label{remark}
As mentioned in the proof of Proposition~\ref{prop:estrem}, for any finite time $t$, including $t$ in the constant $c_\beta$, one may get rid of the power~$1-\alpha$ in the previous proposition, yielding a convergence rate in $\varepsilon^{-c_\beta}$.
\end{rema}

\vspace{0cm}
\bibliographystyle{abbrv}  
\bibliography{arefs}



\end{document}